\documentclass[envcountsame, final]{llncs}
\usepackage{cite}
\usepackage{amsmath,amssymb,bbm}
\usepackage{amsfonts}
\usepackage{mathrsfs}
\usepackage{mathtools}
\usepackage{framed}
\usepackage{subfig}
\usepackage{graphicx}
\usepackage{multirow}
\graphicspath{{./Figures/}}
\usepackage{microtype}
\usepackage{algorithmic}
\usepackage{algorithm}
\usepackage{hyperref}

\providecommand{\mathup}{\mathrm}

\def\R{\mathbb{R}}
\newcommand{\dd}[2]{\frac{d#1}{d#2}}
\newcommand{\grad}[0]{\nabla}
\newcommand{\Diff}{\operatorname{Diff}}

\newcommand{\Dens}{\operatorname{Dens}}
\newcommand{\vol}{\mathup{vol}}

\newcommand{\GI}{G^I}

\begin{document}

\title{Weighted Diffeomorphic Density Matching with Applications to Thoracic Image Registration}
\author{Caleb Rottman\inst{1}, Martin Bauer\inst{2}, Klas Modin\inst{3}, Sarang Joshi\inst{1}}
\institute{Department of Bioengineering, Scientific Computing and Imaging Institute, University of Utah
\and Fakult\"at f\"ur Mathematik, Universit\"at Wien
\and Department of Mathematical Sciences, Chalmers University of Technology and the University of Gothenburg}

\maketitle

\begin{abstract}
In this article we study the problem of thoracic image registration, in particular the estimation of complex anatomical deformations associated with the breathing cycle.
Using the intimate link between the Riemannian geometry of the space of diffeomorphisms and the space of densities, we develop
an image registration framework that incorporates both the fundamental law of conservation of mass as well as spatially varying tissue compressibility properties.
By exploiting the geometrical structure, the resulting algorithm is computationally efficient, yet widely general. 
\end{abstract}

\section{Introduction}

In this paper we consider the problem of tracking organs undergoing deformations as a result of breathing in the thorax and imaged via computed tomography (CT).
This problem has wide scale medical applications, in particular radiation therapy of the lung where accurate estimation of organ deformations during treatment impacts dose calculation and treatment decisions~\cite{sawant2014investigating,keall2005four,suh2014imrt,geneser2011quantifying}.
The current state-of-the-art radiation treatment planning involves the acquisition of a series of respiratory correlated CT (RCCT) images to build 4D (3 spatial and 1 temporal) treatment planning data sets.
Fundamental to the processing and clinical use of these 4D data sets is the accurate estimation of registration maps that characterize the motion of organs at risk as well as the target tumor volumes.

The 3D image produced from X-ray CT is an image of linear attenuation coefficients.
The linear attenuation coefficient $\mu$ of a material is defined as $\mu = \alpha_m \rho_m$, where $\alpha_m$ is the mass attenuation coefficient of the material and $\rho_m$ is the mass density.
The linear attenuation coefficient is proportional to the true density and therefore exhibits conservation of mass.

Currently, the application of diffeomorphisms in medical image registration is mostly limited to the $L^2$ image action of the diffeomorphism group, which is not a mass-preserving transformation.
Furthermore, the diffeomorphisms estimated from typical image registrations algorithms (such as LDDMM~\cite{Beg2005} or ANTS~\cite{Avants2011}) do not accurately model the varying compressibility of different tissues types.
In thoracic datasets, the lungs are highly compressible.
Conversely, the bronchial tubes and the tissue surrounding the lungs are incompressible.
During inhale, as air enters, the lung volume increases and the lung density decreases, while during exhale lung volume decreases and the lung density increases.
But in both inhale and exhale, the lung mass is conserved.

\begin{figure}
  \centering
  \includegraphics[width=0.32\linewidth]{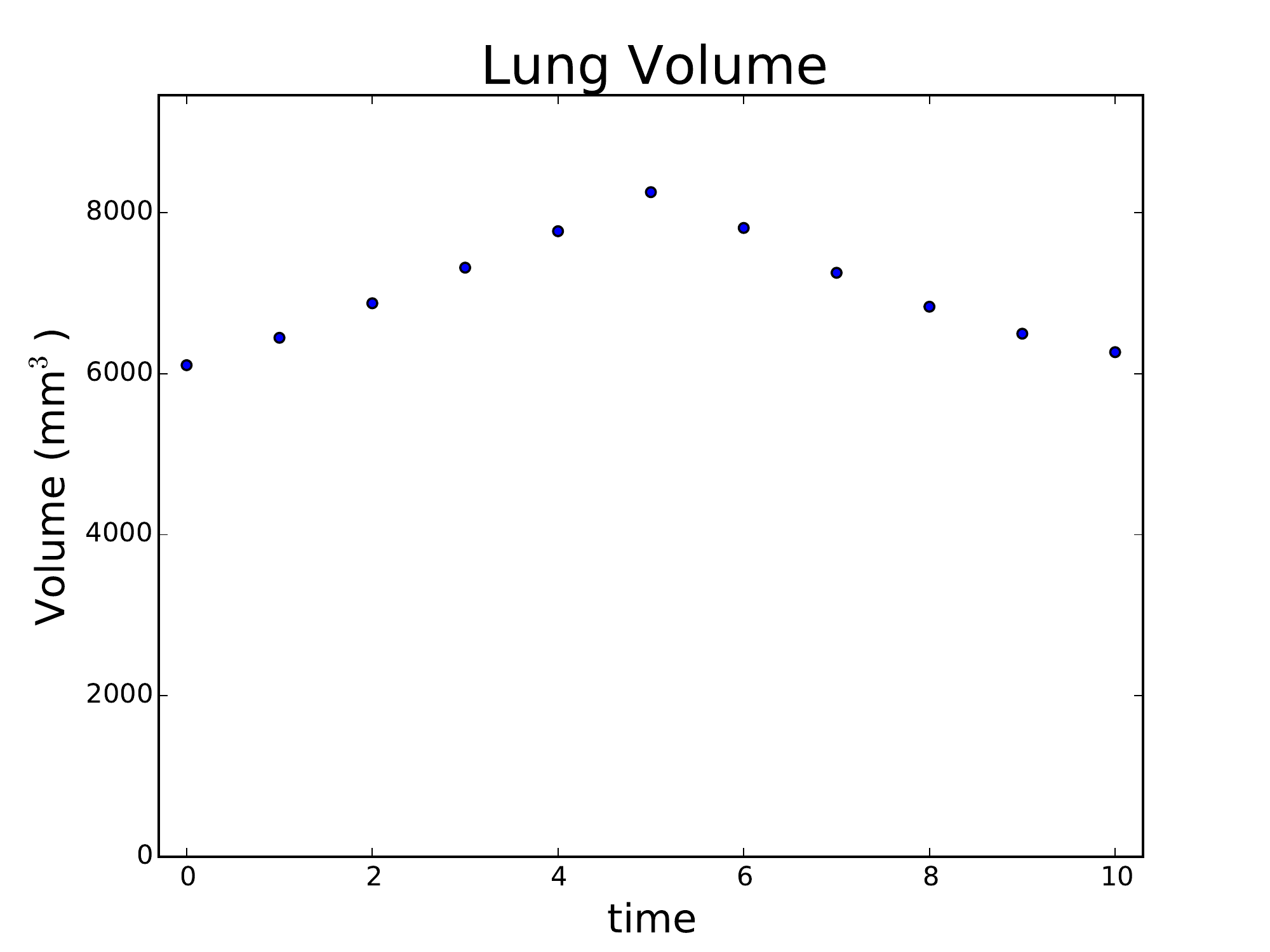}
  \includegraphics[width=0.32\linewidth]{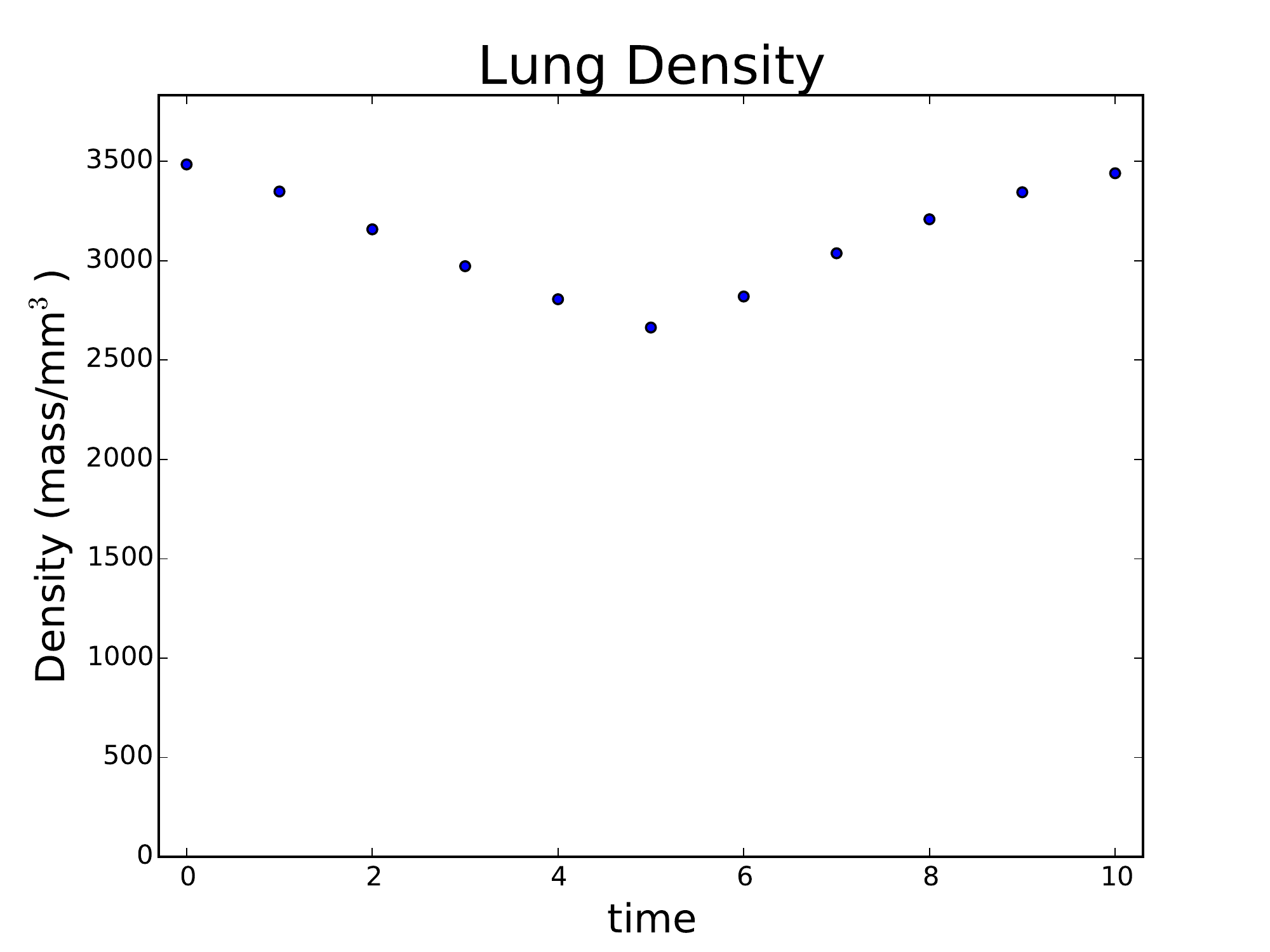}
  \includegraphics[width=0.32\linewidth]{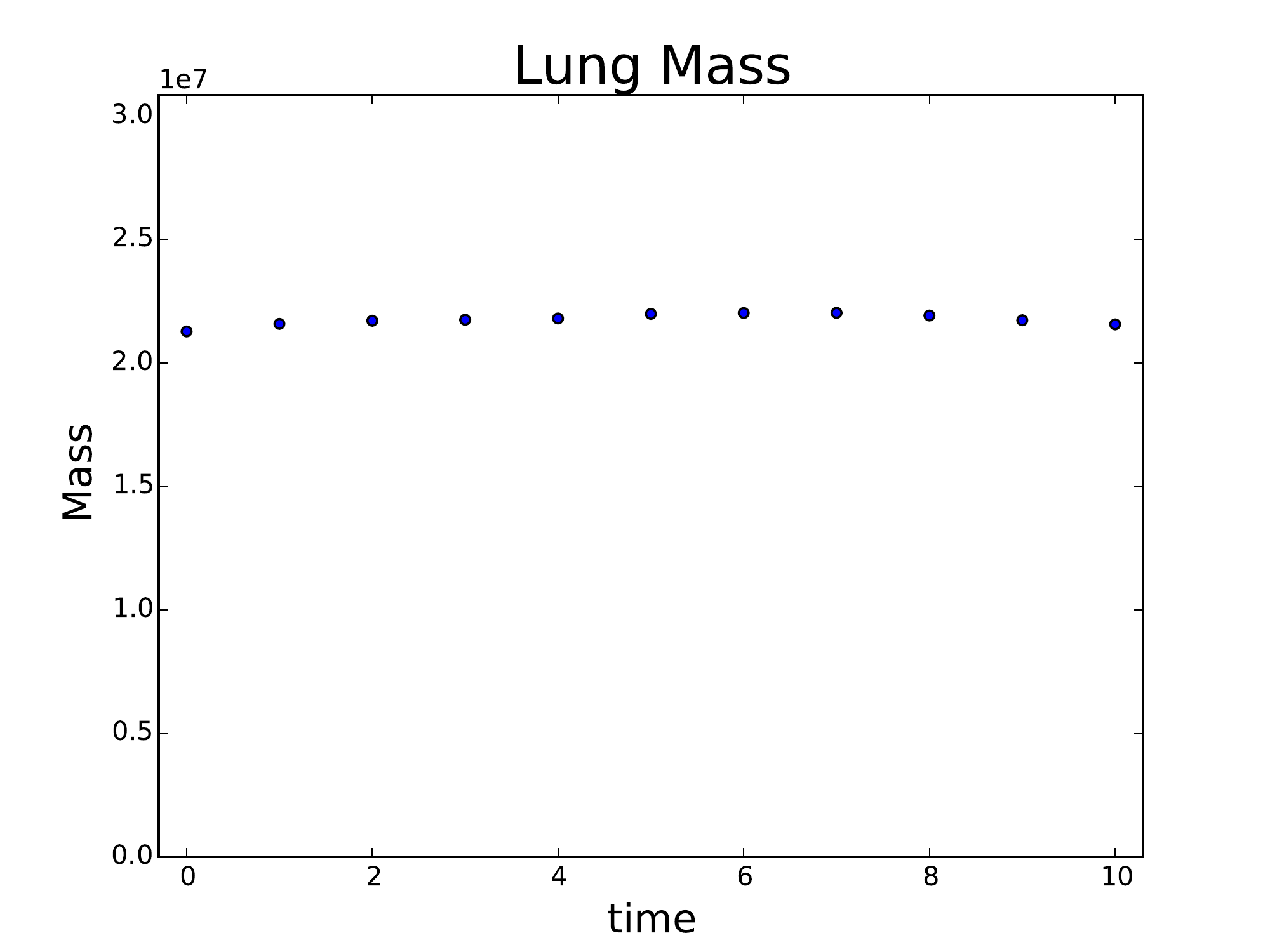}
  \caption{Rat lung data: volume, density, and mass of the lungs during an inhale-exhale breathing cycle.
    As the volume increases, the density decreases, but mass is conserved.}
  \label{fig:MDV}
\end{figure}

In this paper we use a cone-beam CT dataset of a rat acquired at 11 time points of an inhale-exhale breathing cycle.
Figure \ref{fig:MDV} shows the mass, volume, and density of the lungs of a rat at each time point of its breathing cycle, exemplifying these properties.

\begin{figure}
  \centering
  \includegraphics[width=0.55\linewidth]{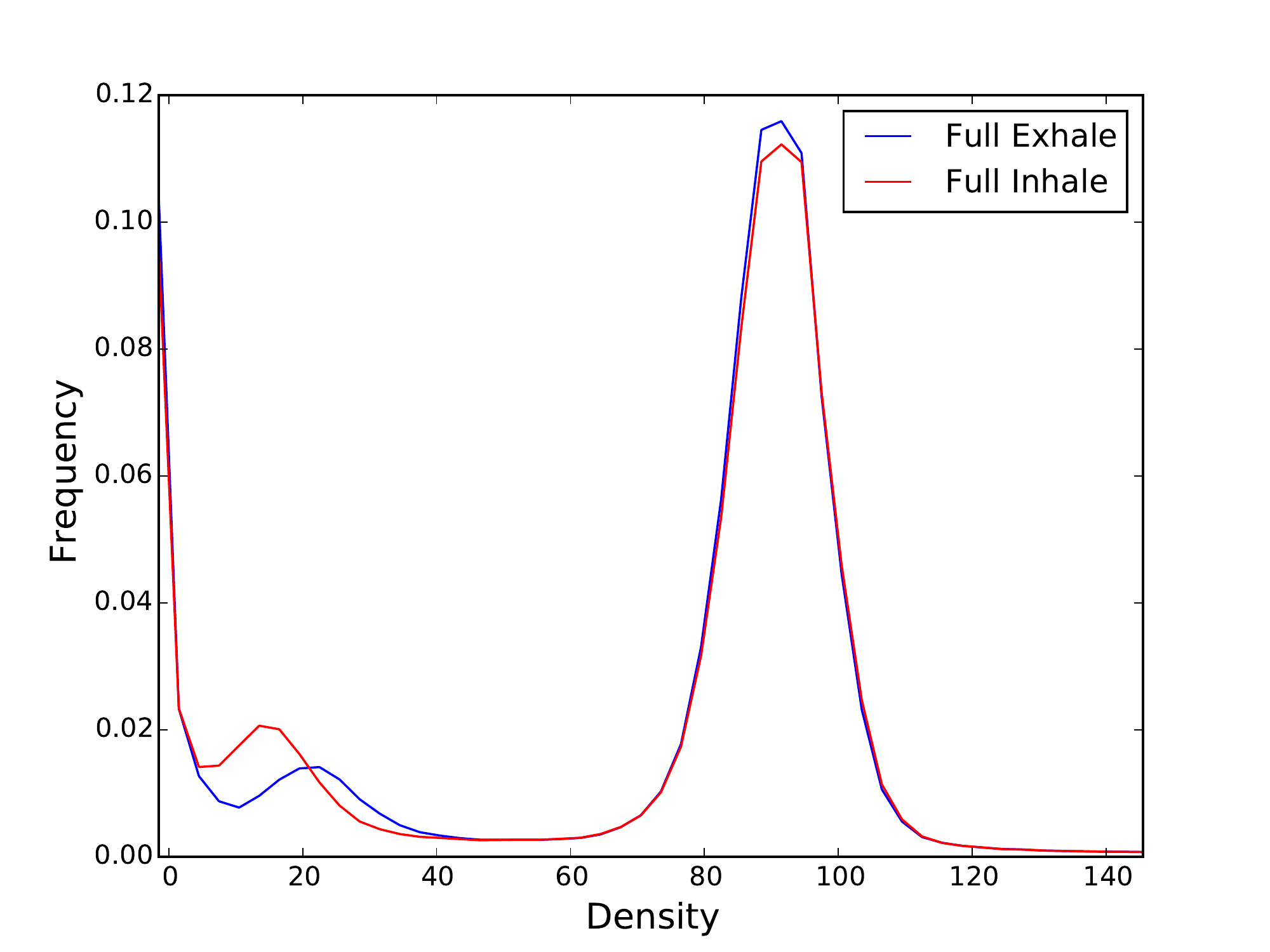}
  \caption{Histograms of a full-inhale and full-exhale image.
    Each histogram has three peaks: the peak at 0 represents surrounding air, the middle peak represents lung tissue, and the peak at 90 represents soft tissue.
    For the lung tissue, the full inhale has higher volume but a lower image intensity than the full exhale, therefore showing conservation of mass.
    For the soft tissue, the average intensity does not change because it is incompressible.
    The slight drop in frequency of the full inhale is due to soft tissue leaving the image boundary.
    }
  \label{fig:hist}
\end{figure}

Both of these effects can be clearly seen in the histograms of a full-inhale and a full-exhale image, as shown in Figure \ref{fig:hist}.

In 2010, the EMPIRE10~\cite{EMPIRE10} challenge compared registration algorithms applied to intra-patient thoracic CT images.
The winner of the competition used an LDDMM method using normalized cross correlation metric \cite{Song2010}.
This method does not model conservation of mass or spatially varying tissue compressibility.
While others in this competition used the density action on these images \cite{Gorbunova2010, Cao2010}, none of these methods incorporate the spatially varying nature of tissue compressibility.

We present an image registration technique that incorporates conservation of mass and organ compressibility.
Instead of the $L^2$ image action of diffeomorphisms, we use the physiologically appropriate density action.
We also regularize the diffeomorphism by using a space-varying penalty which allows for high compressibility of the lung tissue while at the same time enforcing incompressibility of high density structures such as bone.
The algorithm is based on the intimate link between the Riemannian geometry of the space of diffeomorphisms and the space of densities~\cite{KhLeMiPr2013,Mo2015,BaJoMo2015_preprint}.
The resulting algorithm also has the added advantage that it is computationally efficient: orders of magnitude faster than existing diffeomorphic image registration algorithms.

\section{Mathematical Formulation}
Mathematically, the problem is to find a diffeomorphic (bijective and smooth) transformation between two densities on a subset $\Omega\subset\mathbb R^3$.
With a `density' we mean a volume form on $\Omega$, i.e., an element of the form $I\, dx$ where $dx=dx^1\wedge dx^2\wedge dx^3$ is the standard volume element on $\mathbb R^{3}$ and $I=I(x)$ is a non-negative function on $\Omega$.
The space of all densities on $\Omega$ is denoted $\operatorname{Dens}(\Omega)$.
One might, of course, identify $I\,dx$ with its function $I$, and thereby think of $\Dens(\Omega)$ as the set of non-negative functions on $\Omega$.
However, the invariance properties and geometry of the problem are remarkably more transparent when viewing $\Dens(M)$ as a space of volume forms.

The group of diffeomorphisms $\Diff(\Omega)$ acts from the right on $\Dens(\Omega)$ by \emph{pullback:}
the action of $\varphi\in\Diff(\Omega)$ on $I\,d x\in\Dens(\Omega)$ is given by
\begin{equation}
(\varphi,I\, dx)\mapsto \varphi^*(I\,dx) = \left(|D\varphi|\,I\circ\varphi\right)dx,
\end{equation}
where $|D\varphi|$ denotes the Jacobian determinant of $\varphi$.
The corresponding left action is given by \emph{pushforward:}
\begin{equation}
	(\varphi,I\, dx)\mapsto \varphi_*(I\,dx) = (\varphi^{-1})^*(I\,dx) = \left(|D\varphi^{-1}|\,I\circ\varphi^{-1}\right) dx.
\end{equation}

The Riemannian geometry of the group of diffeomorphisms endowed with a suitable Sobolev $H^{1}$ metric is intimately linked to the Riemannian geometry of the space densities with the Fisher--Rao metric.
This has been developed and extensively studied in~\cite{KhLeMiPr2013,Mo2015,BaJoMo2015_preprint}: the basic observation is that there are Sobolev $H^{1}$-metrics on the space of diffeomorphisms that descend to the Fisher--Rao metric on the space of densities.

The distance associated with the Fisher--Rao metric is traditionally defined between \emph{probability densities} (densities of total mass $1$) and is given by
\begin{equation}\label{eq:FR_distance_finite}
	d_F(\mu_0,\mu_1)= \sqrt{\vol(\Omega)}\arccos\left(\frac{1}{\vol(\Omega)}\int_\Omega \sqrt{\frac{\mu_0}{d x}\frac{\mu_1}{d x}} dx\right) \ ,
\end{equation}
where $\mu_0$ and $\mu_1$ are probability densities.
It naturally extends to the space of all densities and the case when $\vol(\Omega)=\infty$, for which it is given by
\begin{equation}\label{eq:FR_distance_infinite}
	d_F^{2}(I_0\,dx,I_1\, dx)= \int_{\Omega} (\sqrt{I_0}-\sqrt{I_1})^{2} dx \, .
\end{equation}
Notice that $d_F^{2}(\cdot,\cdot)$ in this case is the \emph{Hellinger distance}.
For details, see~\cite{BaJoMo2015_preprint}.

The Fisher--Rao metric is the unique Riemannian metric on the space of probability densities that is invariant under the action of the diffeomorphism group~\cite{BBM2015,AJLS2014}.
This invariance property extends to the induced distance function, so
\begin{equation}\label{eq:invariance}
	d_F^{2}(I_0\,dx,I_1\,dx)= d_F^{2}(\varphi_*(I_0\,dx),\varphi_*(I_1\,dx)) \qquad \forall \varphi \in \Diff(\Omega)\;.
\end{equation}

\begin{figure}[htp]
	\centering
	\includegraphics[width=0.7\textwidth]{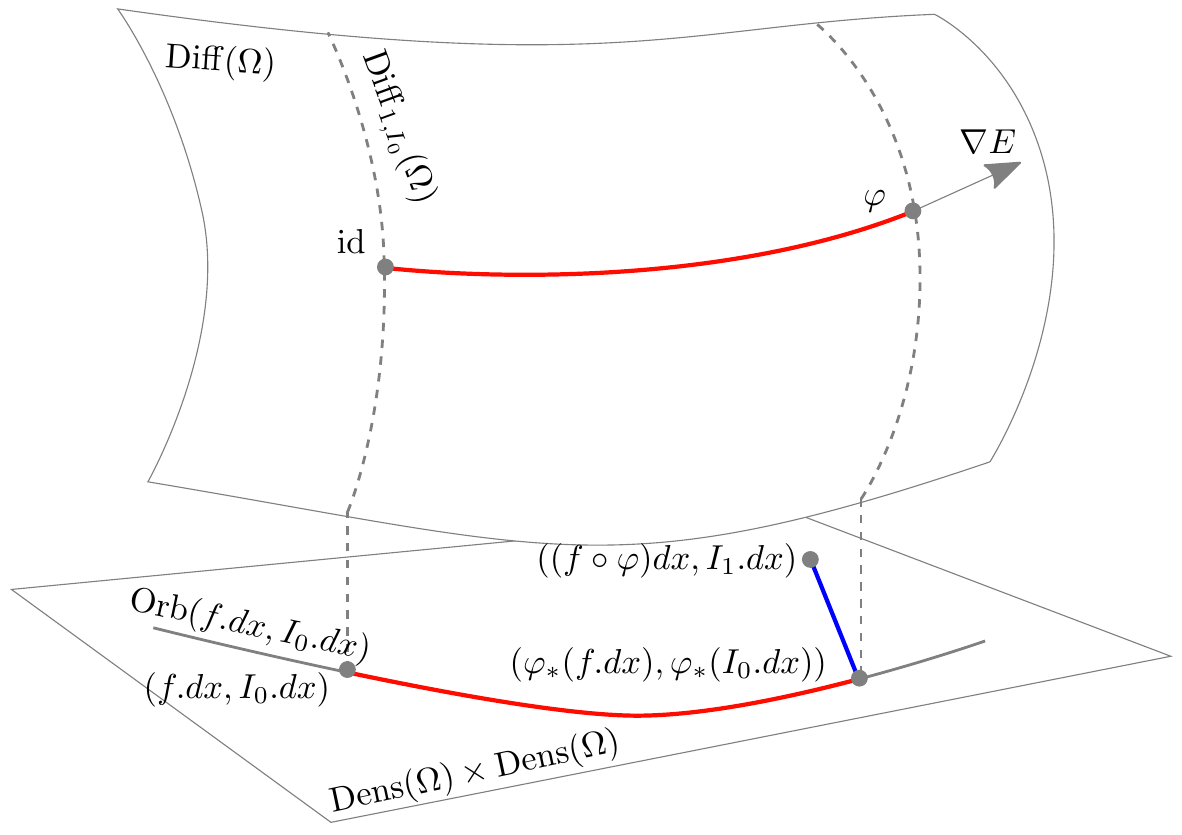}
	\caption{Illustration of the geometry associated with the density matching problem.
	The gradient flow on $\operatorname{Diff}(\Omega)$ descends to a gradient flow on the orbit $\mathrm{Orb}(f\, dx,I_0 \, dx)$.
	While constrained to $\mathrm{Orb}(f\,dx,I_0\,dx)\subset \operatorname{Dens}(\Omega)\times\operatorname{Dens}(\Omega)$, this flow strives to minimize the product Fisher-Rao distance to $((f\circ\varphi)\,dx,I_1\,dx)$.
	}
	\label{Fig1}
\end{figure}

Motivated by the aforementioned properties, we develop a weighted diffeomorphic matching algorithm for matching two density images.
The algorithm is based on the Sobolev $H^{1}$ gradient flow on the space of diffeomorphisms that minimizes the energy functional
\begin{equation}\label{eq:E}
E(\varphi) = d^{2}_F( \varphi_*(f\, dx), (f \circ \varphi^{-1}) dx) + d^{2}_F( \varphi_*(I_0\, dx), I_1\, dx)) .
\end{equation}
This energy functional is only a slight modification of the energy functional studied in~\cite{BaJoMo2015_preprint}.
Indeed, if $f$ in the above equation is a constant $\sigma>0$, then~\eqref{eq:E} reduces to the energy functional of Bauer, Joshi, and Modin~\cite[\S\!~5.1]{BaJoMo2015_preprint}.
Moreover, the geometry described in~\cite[\S\!~5.3]{BaJoMo2015_preprint} is valid also for the functional~\eqref{eq:E}, and, consequently, the algorithm developed in~\cite[\S\!~5.2]{BaJoMo2015_preprint} can be used also for minimizing~\eqref{eq:E}.
There the authors view the energy functional as a constrained minimization problem on the product space $\Dens(\Omega)\times\Dens(\Omega)$ equipped with the
product distance, cf. Fig~\ref{Fig1} and ~\cite[\S\!~5]{BaJoMo2015_preprint} for details on the resulting geometric picture. Related work on  diffeomorphic  density matching using the Fisher Rao metric can be found in \cite{STFV2013,SHV2013}.

Using the invariance property of the Fisher-Rao metric and assuming infinite volume, the main optimization problem associated with the energy functional~\eqref{eq:E} is the following.
\begin{framed}
	\noindent Given densities $I_0\,dx$, $I_1\,dx$, and $f\, dx$, find $\varphi\in\operatorname{Diff}(\Omega)$ minimizing
	\begin{equation}\label{eq:main_problem}
		E(\varphi) = \underbrace{\int_{\Omega} (\sqrt{|D\varphi^{-1}|}-1)^2\,f\circ\varphi^{-1} \,dx}_{E_1(\varphi)} +
		\underbrace{\int_{\Omega}\Big(\sqrt{|D\varphi^{-1}|I_0\circ\varphi^{-1}}- \sqrt{I_1} \Big)^2\,dx}_{E_2(\varphi)} \;.
	\end{equation}
\end{framed}
The invariance of the Fisher-Rao distance can be seen with a simple change of variables $x \mapsto \varphi(y)$, $dx \mapsto |D\varphi| dy$, and $|D\varphi^{-1}| \mapsto \frac{1}{|D\varphi|}$.
Then, Equation \ref{eq:main_problem} becomes
\begin{align}
  \label{eq:main_problemCV}
 E(\varphi) = \int_{\Omega} (1-\sqrt{|D\varphi|})^2\,f \,dy +
 \int_{\Omega}\Big(\sqrt{I_0}- \sqrt{|D\varphi|I_1\circ\varphi} \Big)^2\,dy \;.
\end{align}
To better understand the energy functional $E(\varphi)$ we consider the two terms separately.
The first term $E_1(\varphi)$ is a \emph{regularity measure} for the transformation.
It penalizes the deviation of the diffeomorphism $\varphi$ from being volume preserving.
The density $f\,dx$ acts as a weighting on the domain $\Omega$.
That is, change of volume (compression and expansion of the transformation $\varphi$) is penalized more in regions of $\Omega$ where $f$ is large.
The second term $E_2(\varphi)$  penalizes \emph{dissimilarity} between $I_0\, dx$ and $\varphi^*(I_1\, dx)$.
It is the Fisher--Rao distance between the initial density $I_0\, dx$ and the transformed target density $\varphi^*(I_1\, dx)$.
Because of the invariance~\eqref{eq:invariance} of the Fisher--Rao metric, this is the same as the Fisher--Rao distance between $I_1\, dx$ and $\varphi_*(I_0\, dx)$.


Solutions to problem~\eqref{eq:main_problem} are \emph{not} unique.
To see this, let $\Diff_{I}(\Omega)$ denote the space of all diffeomorphisms preserving the volume form $I\,dx$:
\begin{equation}
	\Diff_{I}(\Omega) = \{ \varphi\in\Diff(\Omega)\mid \; |D\varphi|\,(I\circ\varphi)=I \}	.
\end{equation}
If $\varphi$ is a minimizer of $E(\cdot)$, then $\psi\circ\varphi$ for any
\begin{equation}
	\psi\in\Diff_{1,I_0}(\Omega) \coloneqq  \Diff_{1}(\Omega)\cap \Diff_{I_0}(\Omega)
\end{equation}
is also a minimizer.
Notice that this space is not trivial.
For example, any diffeomorphism generated by a \emph{Nambu--Poisson vector field} (see~\cite{Na1999}), with $I_0$ as one of its Hamiltonians, will belong to it.
A strategy to handle the degeneracy was developed in~\cite[\S\!~5]{BaJoMo2015_preprint}: the fact that the metric is descending with respect to the $H^1$ metric on $\Diff(\Omega)$ can be used to ensure that the gradient flow is \emph{infinitesimally optimal}, i.e., always orthogonal to the null-space.
We employ the same strategy in this paper.
The corresponding geometric picture can be seen in Fig.~\ref{Fig1}.

\section{Gradient Flow Algorithm Development}
We now derive in detail the algorithm used to optimize the functional defined in Equation~\ref{eq:main_problemCV}.
The $H^{1}$-metric on the space of diffeomorphisms is defined using the Hodge laplacian on vector fields and is given by:
\begin{equation}\label{eq:H1dot}
	\GI_{\varphi}(U,V) = \int_{\Omega}\langle -\Delta u,v\rangle dx \ .
\end{equation}
Due to its connections to information geometry we also refer to this metric as \emph{information metric}.
 Let $\nabla^{\GI} E$ denote the gradient with respect to the information metric defined above. Our approach to minimize the functional of~\eqref{eq:main_problemCV} is to use a simple Euler integration of the discretization of the gradient flow:
 \begin{equation}
 {\dot \varphi } = - \nabla^{\GI} E(\varphi)
 \end{equation}
 The resulting final algorithm (Algorithm \ref{algorithm}) is order of magnitudes faster than LDDMM, since we are not required to time integrate the geodesic equations, as necessary in LDDMM~\cite{Younes2009S40}.

 In the following theorem we calculate the gradient of the energy functional:
 \begin{theorem}
  The $\GI$--gradient of the matching functional~\eqref{eq:main_problemCV} is given by
  \begin{align}
  \nabla^{\GI} E &= -\Delta^{-1} \Big( -\grad\big(f\circ \varphi^{-1} (1-\sqrt{|D\varphi^{-1}|})\big) - \nonumber \\
  & \hspace{30pt} \sqrt{|D\varphi^{-1}|\,I_0\circ\varphi^{-1}} \grad\big(\sqrt{I_1}\big) + \grad\big(\sqrt{|D\varphi^{-1}|\,I_0\circ\varphi^{-1}}\big)\sqrt{I_1} \Big)\;.
\end{align}
 \end{theorem}
 \begin{remark}
Notice that in the formula for $\nabla^{\GI} E$ we never need to compute $\varphi$, so in practice we only compute $\varphi^{-1}$.
We update this directly via $\varphi^{-1}(y) \mapsto \varphi^{-1}(y + \epsilon \nabla^{\GI} E)$ for some step size $\epsilon$.
\end{remark}
\begin{proof}
We first calculate the variation of the energy functional. Therefore let $\varphi_s$ be a family of diffeomorphisms parameterized by the real variable $s$, such that
\begin{equation}
  \varphi_0 = \varphi \hspace{10pt} \mathrm{and} \hspace{10pt} \dd{}{s} \Big|_{s=0} \varphi_s = v \circ \varphi.
\end{equation}
We use the following identity, as derived in \cite{hinkle2013idiff}:
\begin{align}
  \dd{}{s} \Big|_{s=0} \sqrt{|D\varphi_s|} =& \frac{1}{2}\sqrt{|D\varphi|}\mathrm{div}(v) \circ \varphi .
\end{align}
The variation of the first term of the energy functional is
\begin{align}
    \dd{}{s} \Big|_{s=0} E_1(\varphi) =& \int_\Omega f(x) (\sqrt{|D\varphi(x)|}-1)\sqrt{|D\varphi(x)|} \mathrm{div}(v) \circ \varphi(x) dx
\end{align}
We do a change of variable $x \mapsto \varphi^{-1}(y)$, $dx \mapsto |D\varphi^{-1}(y)| dy$, using the fact that $|D\varphi(x)| = \frac{1}{|D\varphi^{-1}(y)|}$;
\begin{align}
  =& \int_\Omega f\circ \varphi^{-1}(y) (1-\sqrt{|D\varphi^{-1}(y)|}) \mathrm{div}(v)(y) dy \\
  =& \left \langle f\circ \varphi^{-1} (1-\sqrt{|D\varphi^{-1}|}), \mathrm{div}(v) \right \rangle_{L^2(\R^3)} \\
  =& -\left \langle \grad\left(f\circ \varphi^{-1} (1-\sqrt{|D\varphi^{-1}|})\right), v \right \rangle_{L^2(\R^3)}
\end{align}
using the fact that the adjoint of the divergence is the negative gradient.
For the second term of the energy functional, we expand the square
\begin{align}
  E_2(\varphi) &= \int_\Omega  I_0(x) - 2\sqrt{I_0(x) I_1 \circ \varphi(x) |D\varphi(x)|} + I_1 \circ \varphi(x)|D\varphi(x)| dx
\end{align}
Now $\int_\Omega I_1 \circ \varphi(x)|D\varphi(x)| dx$ is constant (conservation of mass), so we only need to minimize over the middle term.
The derivative is then
\begin{align}
  \dd{}{s} \Big|_{s=0} E_2(\varphi) &= - \int_\Omega 2 \sqrt{I_0(x)} \big(\grad \sqrt{I_1}^T v\big)\circ \varphi(x) \sqrt{|D\varphi(x)|} \nonumber \\
 & \hspace{20pt}  - \sqrt{I_0(x)I_1 \circ \varphi(x) |D\varphi(x)|} \mathrm{div}(v) \circ \varphi(x) dx.
\end{align}
We do the same change of variables as before:
\begin{align}
   &= - \int_\Omega \sqrt{I_0\circ\varphi^{-1}(y)} \frac{|D\varphi^{-1}(y)|}{\sqrt{|D\varphi^{-1}(y)}|} \big( 2\grad \sqrt{I_1(y)}^T v(y) +\sqrt{I_1(y)}\mathrm{div}(v)(y) \big) \\
   &= -\left \langle 2 \sqrt{|D\varphi^{-1}|\,I_0\circ\varphi^{-1}} \grad \sqrt{I_1}, v \right \rangle_{L^2(\R^3)} \nonumber \\
   &\hspace{50pt} - \left \langle  \sqrt{|D\varphi^{-1}|\,I_0\circ\varphi^{-1} I_1}, \mathrm{div}(v) \right \rangle_{L^2(\R^3)} \\
   &= \left \langle - \sqrt{|D\varphi^{-1}|\,I_0\circ\varphi^{-1}} \grad \sqrt{I_1}, v \right \rangle_{L^2(\R^3)} \nonumber \\
    &\hspace{50pt} + \left \langle \grad\left(\sqrt{|D\varphi^{-1}|\,I_0\circ\varphi^{-1}}\right)\sqrt{I_1}, v \right \rangle_{L^2(\R^3)}.
\end{align}
From the above equations we conclude that:
\begin{align}
  \label{eq:finalgrad}
  -\Delta(\nabla^{\GI} E) &= -\grad\left(f\circ \varphi^{-1} (1-\sqrt{|D\varphi^{-1}|})\right) \nonumber \\
  & \hspace{10pt} - \sqrt{|D\varphi^{-1}|\,I_0\circ\varphi^{-1}} \grad \sqrt{I_1} + \grad\left(\sqrt{|D\varphi^{-1}|\,I_0\circ\varphi^{-1}}\right)\sqrt{I_1}
\end{align}
Since we are taking the Sobolev gradient of $E$, we apply the inverse Laplacian to the right hand side of Equation \ref{eq:finalgrad} to solve for $\nabla^{\GI} E$.
\end{proof}

\begin{algorithm}
  \caption{Final Algorithm}\label{alg:main}
  \label{algorithm}
  \begin{algorithmic}
    \STATE{Choose $\epsilon > 0$}
    \STATE{Set $\varphi^{-1} = \mathrm{id}$}
    \STATE{Set $|D\varphi^{-1}| = 1 $}
    \FOR{$iter = 1...\mathrm{NumIters}$}
    \STATE{Compute $\varphi_*I_0 = I_0 \circ \varphi^{-1} |D\varphi^{-1}|$}
    \STATE{Compute $u = -\grad\big(f\circ \varphi^{-1} (1-\sqrt{|D\varphi^{-1}|})\big) -
      \sqrt{\varphi_*I_0} \grad \sqrt{I_1} + \grad(\sqrt{\varphi_*I_0})\sqrt{I_1}$}
    \STATE{Compute $v = -\Delta^{-1}(u)$}
    \STATE{Update $\varphi^{-1} \mapsto \varphi^{-1}(y + \epsilon v)$}
    \STATE{Update $|D\varphi^{-1}| \mapsto  |D\varphi^{-1}|\circ \varphi^{-1}\mathrm{e}^{-\epsilon \mathrm{div}(v)}$}
    \ENDFOR
  \end{algorithmic}
\end{algorithm}

\begin{remark}
  Algorithm~\ref{alg:main} constructs the mapping $\varphi^{-1}$ by numerically integrating the vector field $v$.
  Thus, for small enough $\epsilon$, the computed transformation $\varphi^{-1}$ is a diffeomorphism (as is also the case in LDDMM).
\end{remark}

\section{Results}

We applied the proposed method to the previously mentioned rat dataset.
In this dataset, an anesthetized rat was placed on a mechanical ventilator.
This ventilator sent 11 gate signals to the cone-beam CT per breathing cycle, assuring that all projections would all be acquired at a consistent points of the breathing cycle \cite{Jacob2011}.
Previous literature has shown that cone-beam CT is inadequate in estimating the true linear attenuating coefficient density \cite{DeVos2009}, so we empirically estimated the density as the square of the the original data.

\begin{figure}[t]
  \centering
  \begin{tabular}{cccc}
    \includegraphics[width=.23\textwidth]{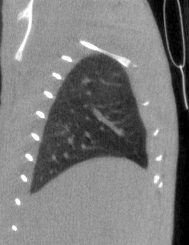}&
    \includegraphics[width=.23\textwidth]{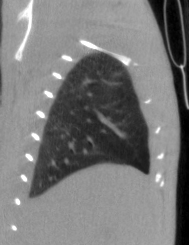}&
    \includegraphics[width=.23\textwidth]{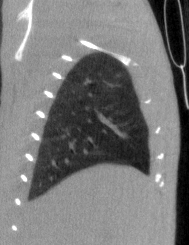}&
    \includegraphics[width=.23\textwidth, trim = 23mm 19mm 47mm 18mm, clip]{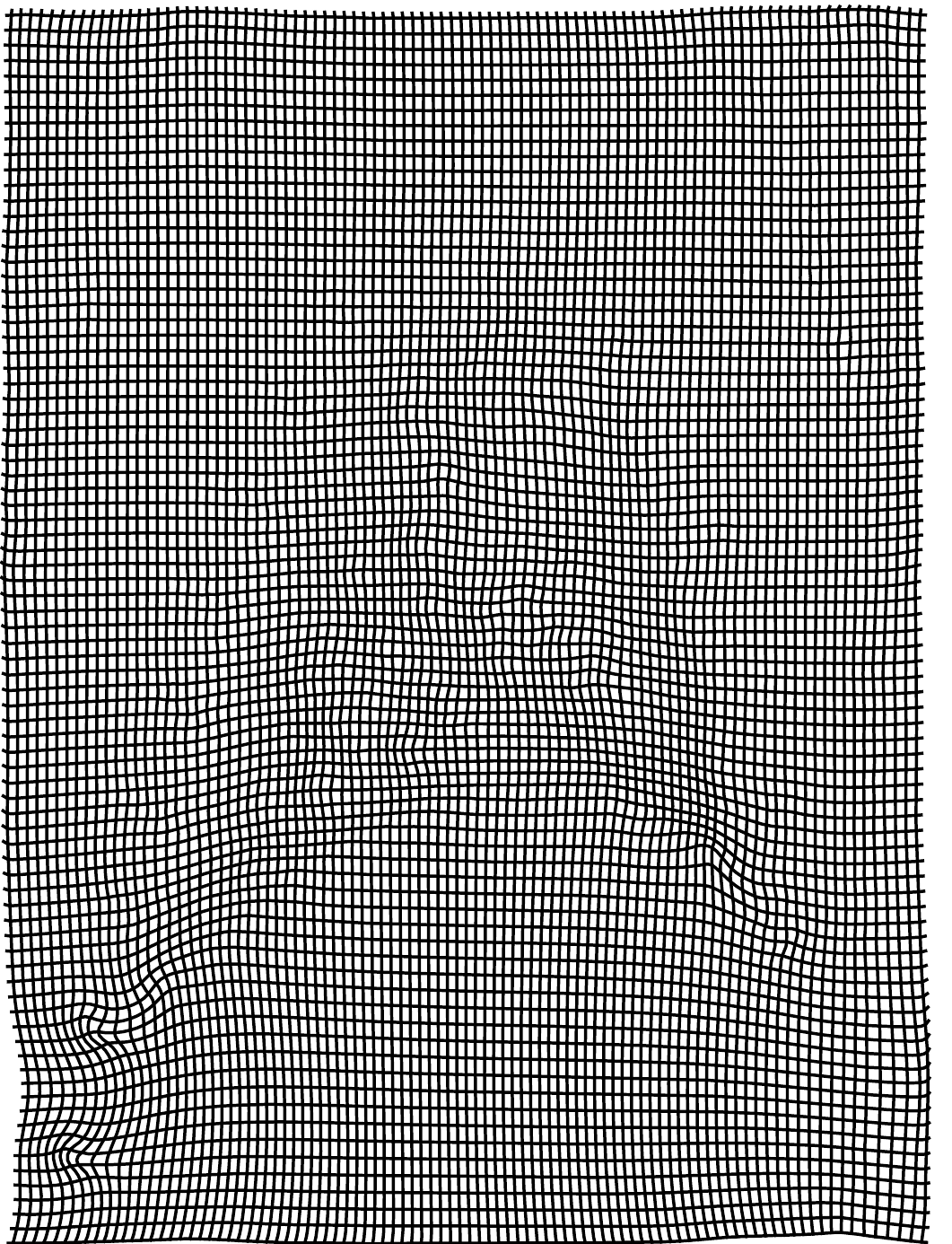}
    \\
    $I_{ex}$ & $\varphi_*(I_{ex}dx)$ & $I_{in}$ & $\varphi^{-1} $
  \end{tabular}
  \caption{Density action results. This figure shows the lung image at the full exhale, the full exhale deformed via the density action, and the corresponding image at full inhale. Shown in the right panel is the estimated deformation. }
  \label{fig:images}
\end{figure}

For these results we estimated the deformation from the full-exhale  to the full-inhale image.
The deformation was computed on the resolution of the original 3D volume ($245\times189\times217$); all the figures show the same 2D coronal slice of this volume.
Shown in Fig.~\ref{fig:images} are the coronal sections of full exhale, the full exhale deformed via the density action, and the corresponding image at full inhale and the estimated deformation.

\begin{figure}[h]
  \centering
  \begin{tabular}{ccc}
    \includegraphics[width=.32\textwidth]{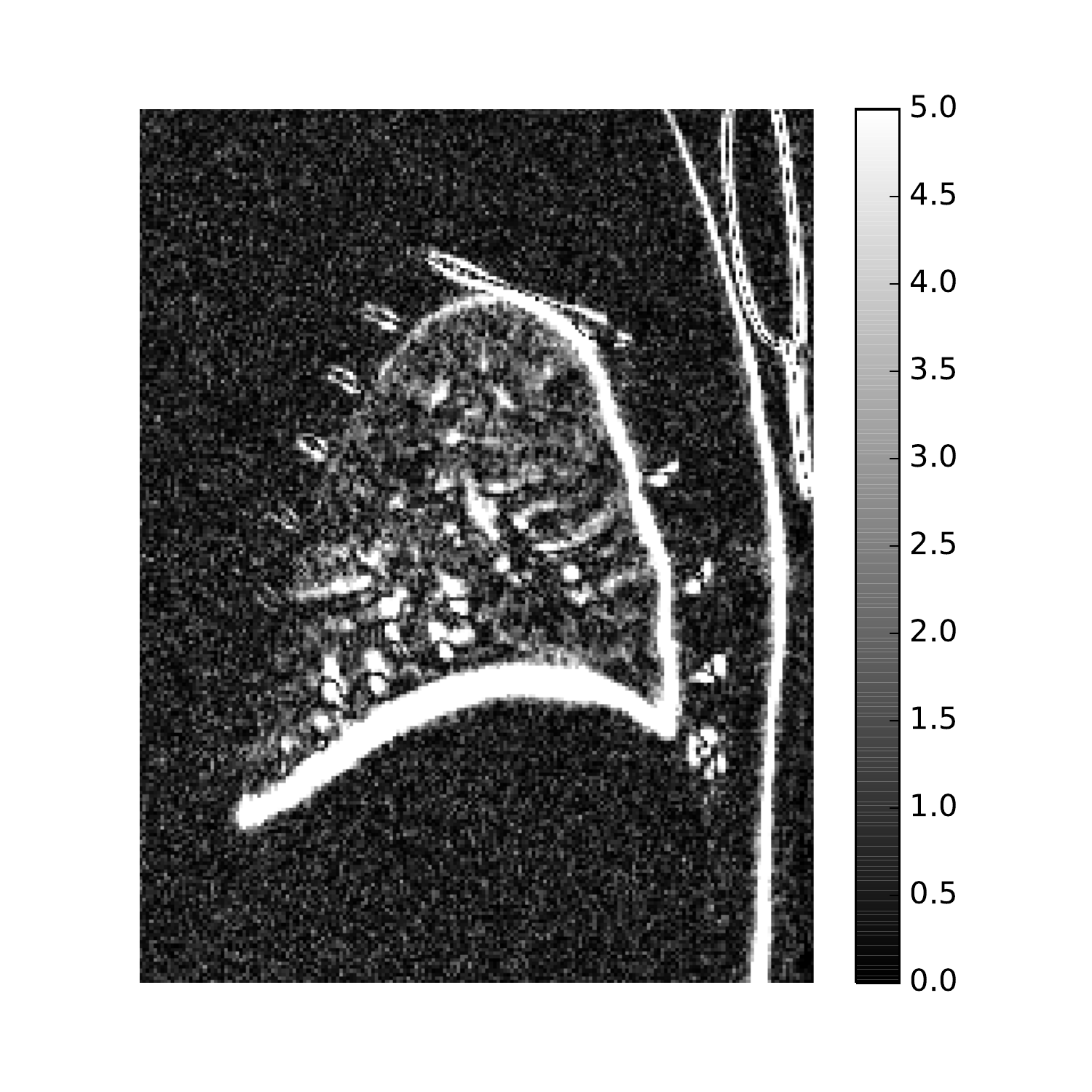} &
    \includegraphics[width=.32\textwidth]{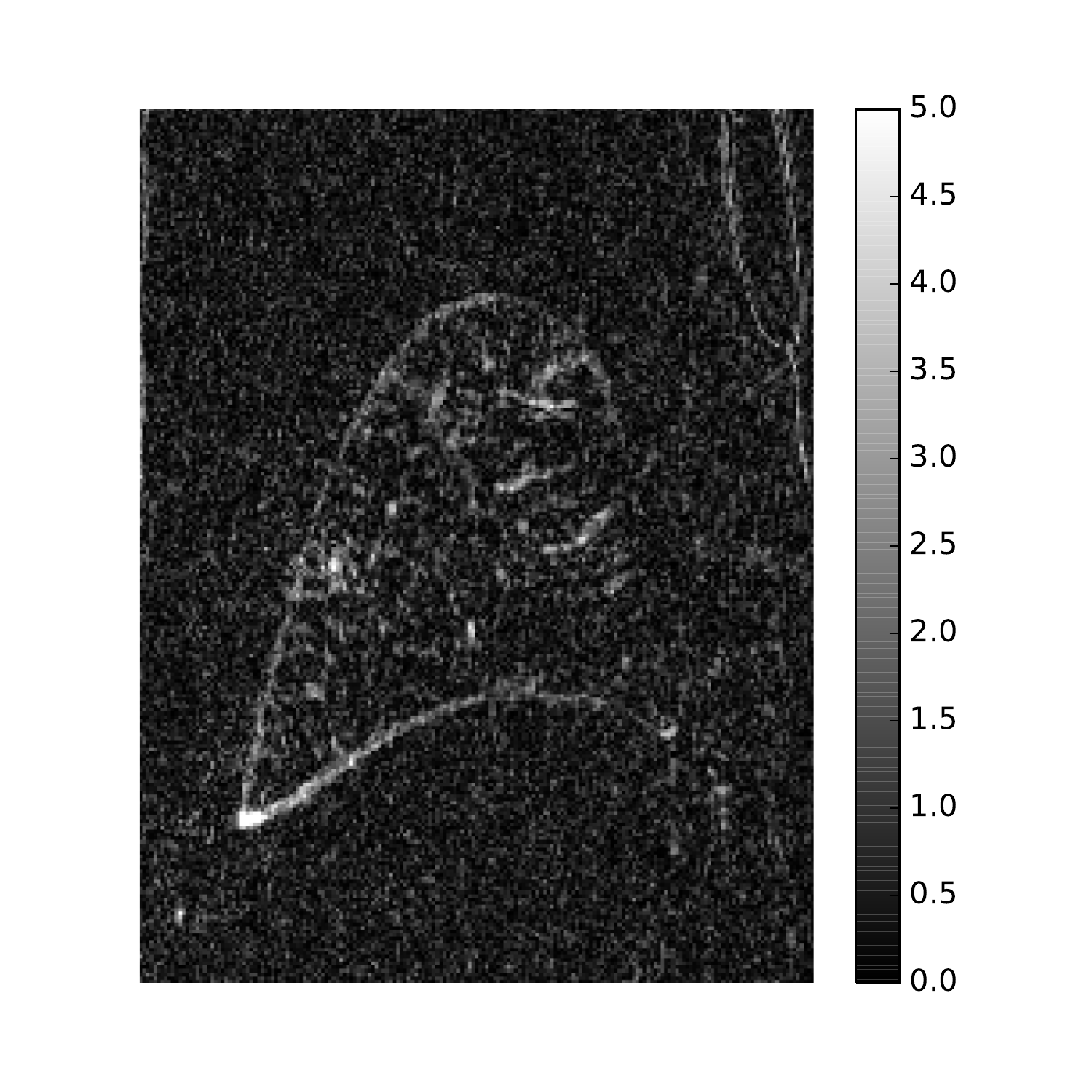} &
    \includegraphics[width=.32\textwidth]{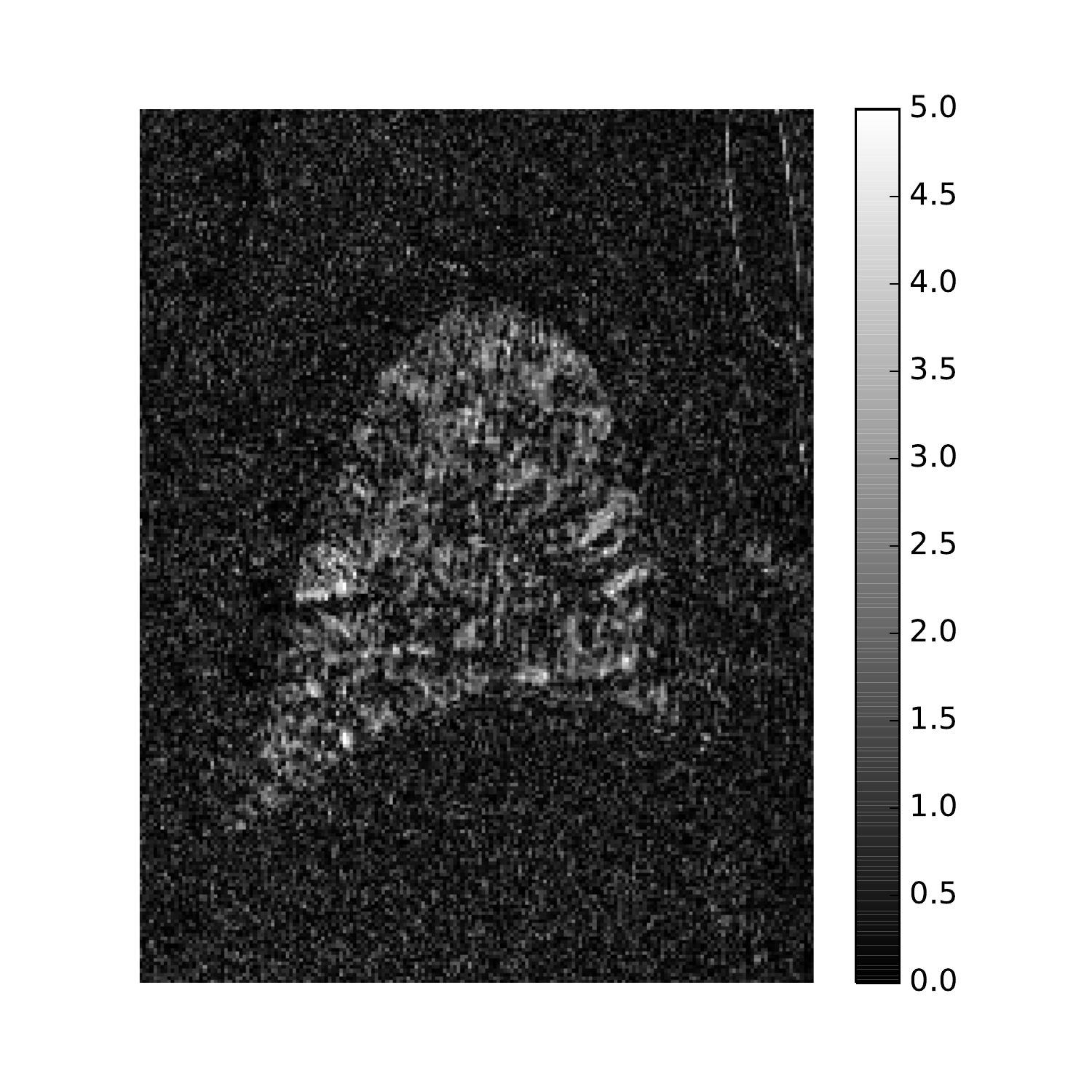} \\
    $|I_{in}-I_{ex}|$ & $|I_{in} -\varphi_*(I_{ex}dx)|$ & $|I_{in} - I_{ex} \circ \varphi^{-1}|$ (LDDMM)
  \end{tabular}
  \caption{Absolute value of image differences:
    The left panel shows the difference between the original full exhale and the full inhale images.
    The center panel shows the result after registration using the proposed method.
    The right image shows the result using LDDMM with the $L^2$ image action.
    In LDDMM, there is significant error inside the lung due to the $L^2$ action not preserving mass.
    }
    \label{fig:diff}
\end{figure}

For the compressibility penalty $f$, we used a soft thresholding of the intensity values of the initial image using the logistic function.
High intensity regions of the CT image (corresponding to bone and soft tissue) were given a high penalty ($f(x) = 10\sigma$) and low intensity regions of the CT image (corresponding to air and lungs) were given a low penalty ($f(x) = .1\sigma)$ (see Figure~\ref{fig:plots})

\begin{figure}[h!]
  \centering
  \begin{tabular}{ccc}
    \includegraphics[width=.32\textwidth]{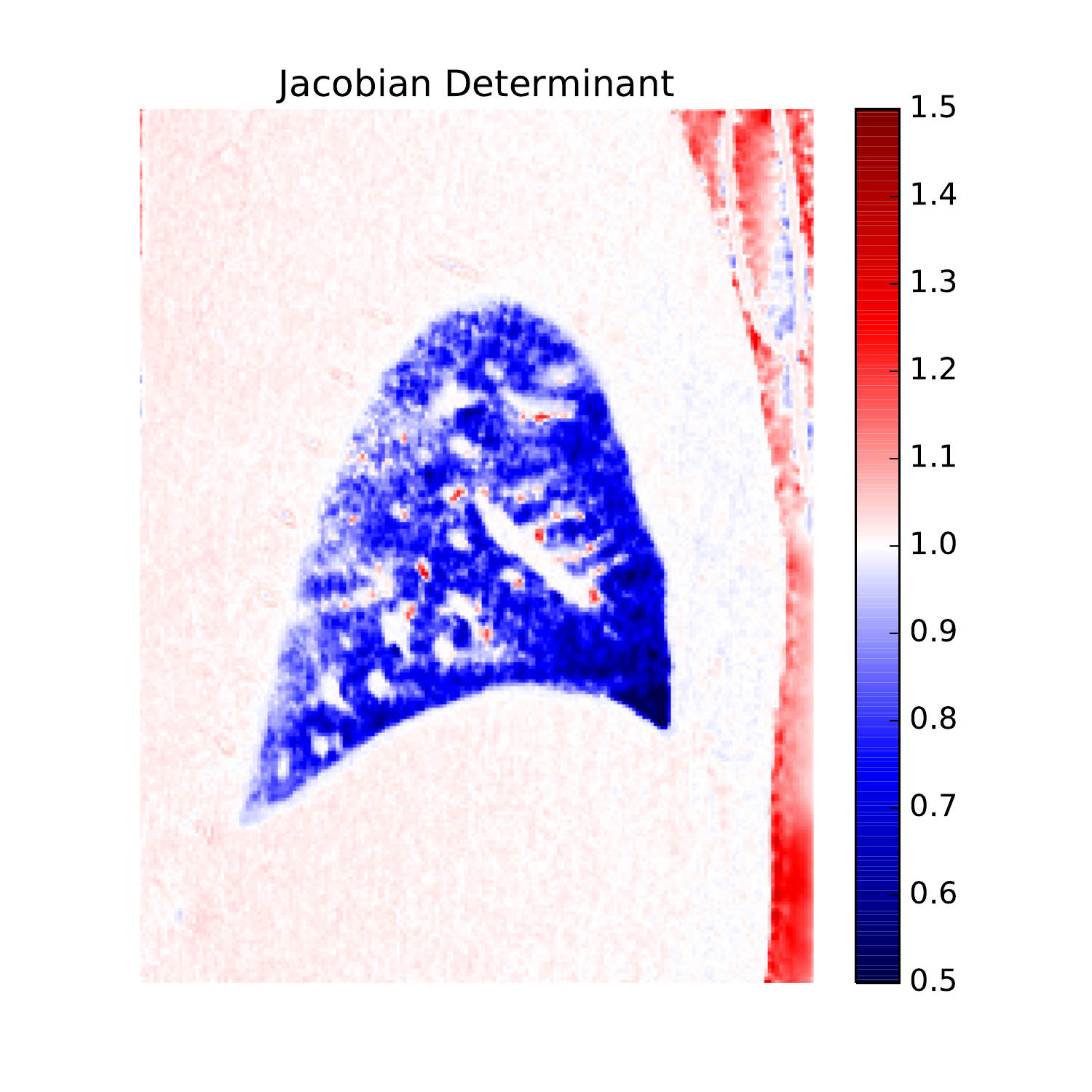} &
    \includegraphics[width=.32\textwidth]{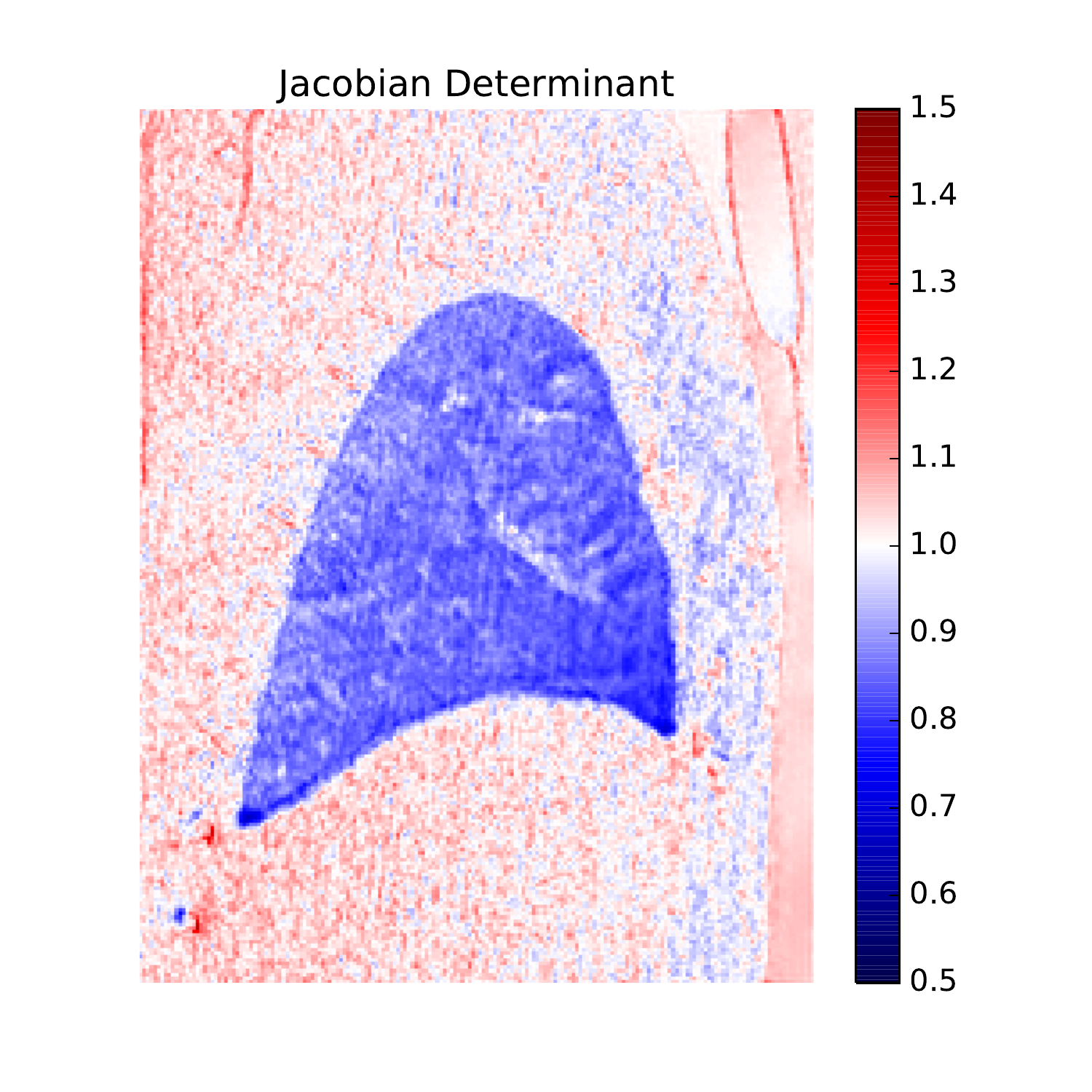} &
    \includegraphics[width=.32\textwidth]{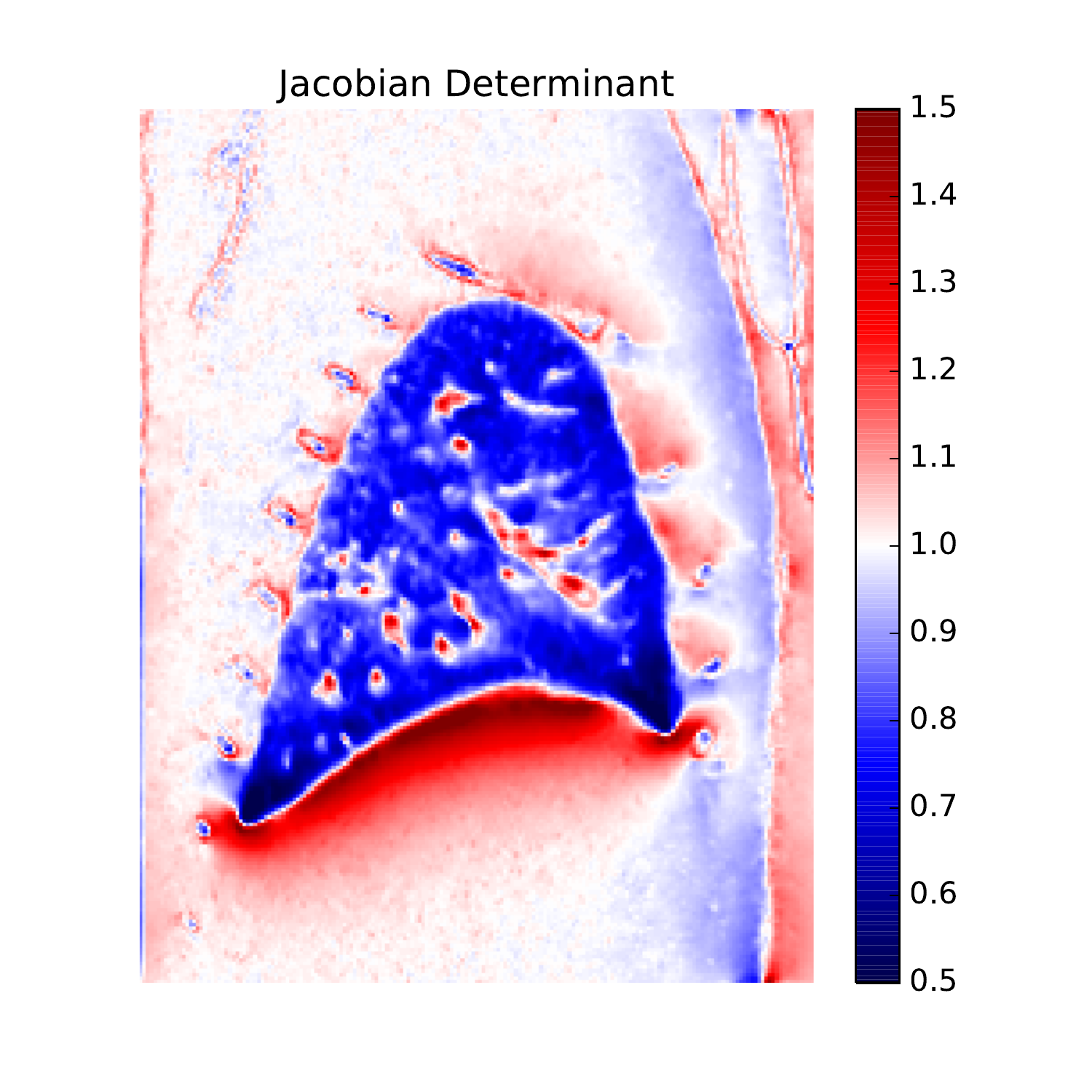} \\
    $f(x) = \mathrm{sig}(I_0(x))$ & $f(x) = 1$ & LDDMM
  \end{tabular}
  \caption{Jacobian determinants:
    On the left is the Jacobian determinant of the transformation estimated by the proposed method.
    Notice that the volume change is confined to inside the lungs and outside the body.
    In the center we use the density action, but without a local-varying penalty (i.e.\, $f(x) = \sigma$).
    On the right is the Jacobian determinant using LDDMM.
    Without the local-varying penalty, there is contraction and expansion outside of the lungs.
    In LDDMM, the contraction and expansion outside of the lungs is even more severe.
  }
  \label{fig:detDphi}
\end{figure}

\begin{figure}[h]
  \centering
  \begin{tabular}{cc}
    \includegraphics[width=.44\textwidth]{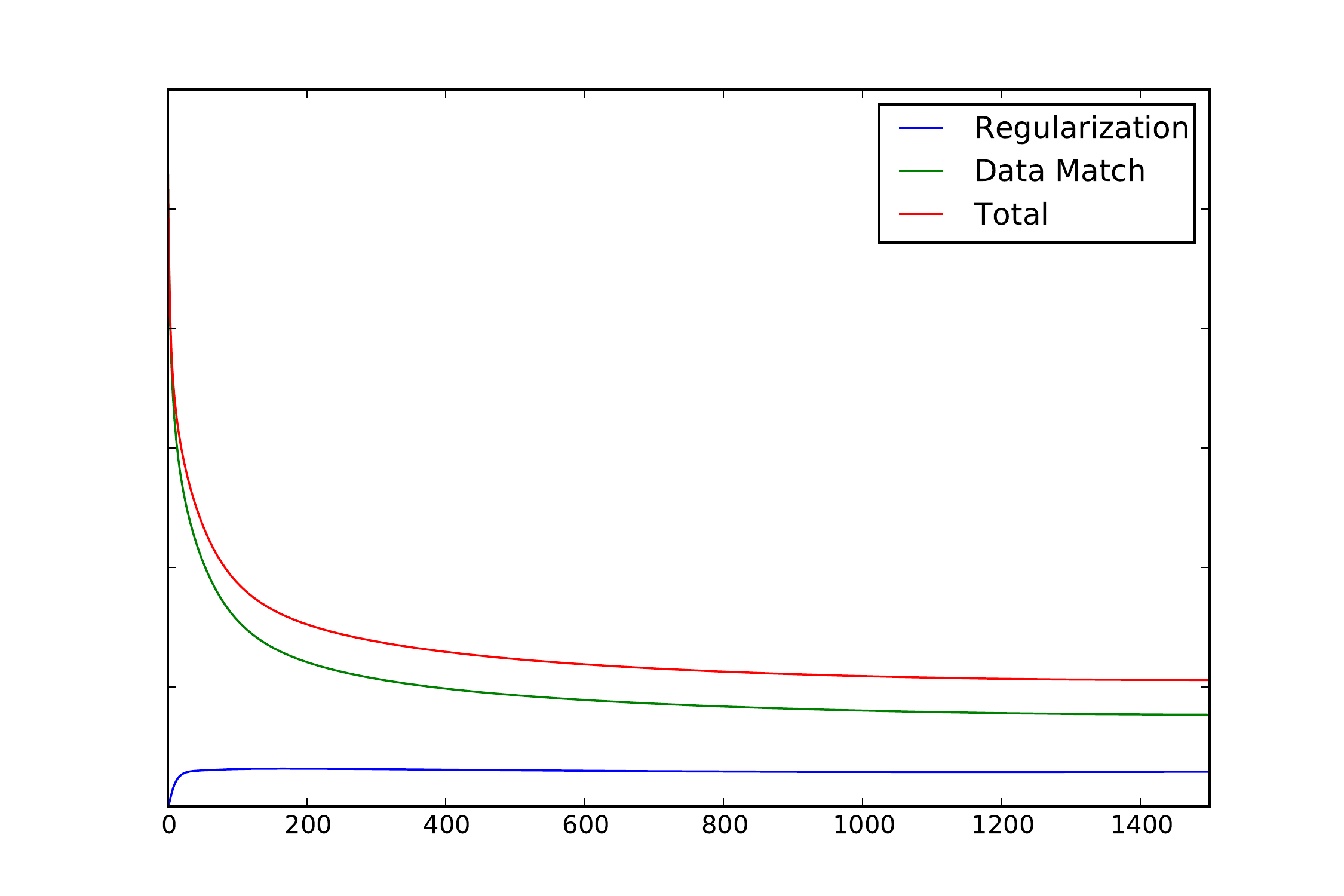} &
    \includegraphics[width=.44\textwidth]{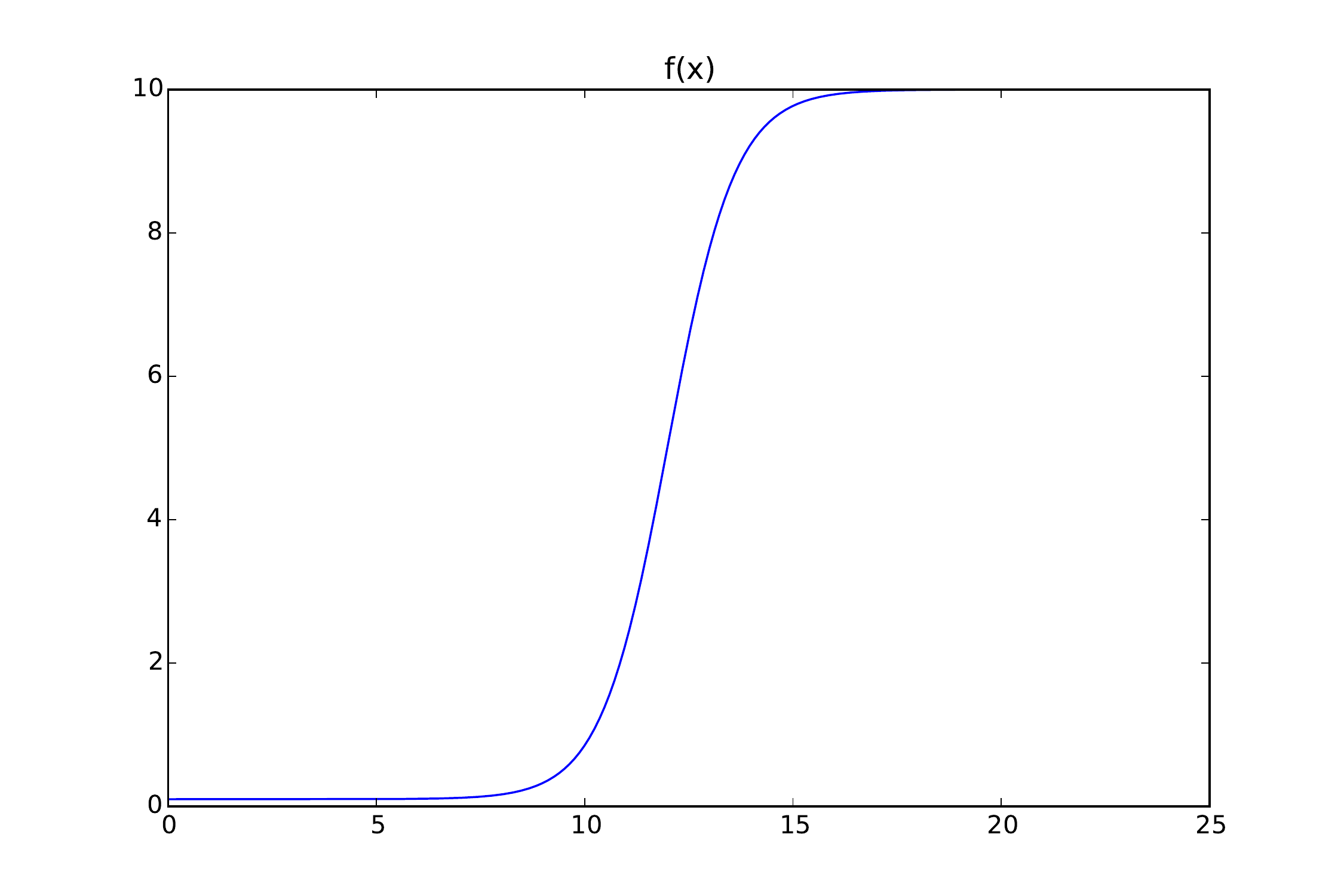} \\
    Energy & $\mathrm{sig}(x)$
  \end{tabular}
  \caption{Energy plot and the logistic function used for the penalty.}
  \label{fig:plots}
\end{figure}

We implemented the proposed algorithm and LDDMM on a single Titan-Z GPU (using the PyCA software package \cite{pyca} \href{http://bitbucket.org/scicompanat/pyca}{bitbucket.org/scicompanat/pyca} ) for comparison.
The difference images are pictured in Figure~\ref{fig:diff}.
The problem of LDDMM using the $L^2$ action can be seen in this image.
The Jacobian determinants are in Figure~\ref{fig:detDphi}.
The proposed method constrains the contraction and expansion to inside the lung and outside the body.
In this figure we also show the results of using the density action with a constant penalty function ($f(x)=\sigma$.).

The proposed algorithm is significantly faster than LDDMM; it runs at 400 iterations per minute while LDDMM runs at 45 iterations per minute.
We used 10 time steps to integrate the geodesic equations associated with the LDDMM formulation.
Since we are not required to integrate the geodesic equations in the proposed algorithm, we have nearly a 10x speedup compared to LDDMM.


\section{Discussion}

In this paper, we introduced a computationally efficient method for estimating registration maps between thoracic CT images.
The proposed solution accurately incorporates the fundamental property of mass conservation and the spatially varying compressibility of thoracic anatomy.
We conserve mass by viewing the images as densities and applying the density action of a diffeomorphism instead of the typical $L^2$ action.
We limit the volume change in incompressible organs by placing a space-varying penalty on the Jacobian determinant of the diffeomorphism.
While any non-negative function $f(x)$ can be used, we simply use a soft-thresholding function on the initial image.
This choice is based on the assumption that low CT image values (such as the lungs and air) exhibit a large amount of volume change whereas high images values (such as other soft tissue and bone) are quite incompressible.

\section*{Acknowledgments}
The authors thank Rick Jacob at the Pacific Northwest National Laboratory for the imaging data which was funded by a grant from the National Heart, Lung, and Blood Institute of the National Institutes of Health (R01 HL073598). 
The works was partially supported by the grant NIH R01 CA169102-01A13, the Swedish Foundation for Strategic Research (ICA12-0052), an EU Horizon 2020 Marie Sklodowska-Curie Individual Fellowship (661482) and by the Erwin Schr\"{o}dinger Institute programme: Infinite-Dimensional Riemannian Geometry with Applications to Image Matching and Shape Analysis. M. Bauer was supported by the European Research Council (ERC), within the project 306445 (Isoperimetric Inequalities and Integral Geometry) and by the FWF-project P24625 (Geometry of Shape spaces).

\bibliographystyle{splncs03}
\bibliography{RBKJ2015}

\end{document}